\documentstyle[12pt,twoside]{amsart}

\title{
A supplement to Fujino's paper: \\
On isolated log canonical singularities with index one}
\author{Shihoko Ishii}
\address{Graduate School of Mathematical Science, University of Tokyo,
3-8-1 Komaba, Meguro, 153-8914 Tokyo, Japan}
\email{{\tt shihoko@@ms.u-tokyo.ac.jp}}

\begin{document}
\maketitle

\newcommand{\bC}{{\Bbb C}}
\newcommand{\bP}{{\Bbb P}}
\newcommand{\bZ}{{\Bbb Z}}
\newcommand{\bQ}{{\Bbb Q}}
\newcommand{\bR}{{\Bbb R}}
\newcommand{\bN}{{\Bbb N}}
\newcommand{\bA}{{\Bbb A}}
\newcommand{\bG}{{\Bbb G}}
\newcommand{\bV}{{\Bbb V}}
\newcommand{\R}{{\cal R}}
\newcommand{\J}{{\cal J}}
\newcommand{\I}{{\cal I}}
\newcommand{\ba}{{\bf a}}
\newcommand{\bp}{{\bf p}}
\newcommand{\bq}{{\bf q}}
\newcommand{\e}{{\cal E}}
\newcommand{\spec}{\operatorname{Spec}}
\newcommand{\Proj}{\operatorname{Proj}}
\newcommand{\Hom}{\operatorname{Hom}}
\newcommand{\ord}{\operatorname{ord}}
\newcommand{\vol}{\operatorname{vol}}
\newcommand{\st}{{\spec k[[t]]}}
\newcommand{\stm}{{\spec k[t]/(t^{m+1})}}
\newcommand{\sTm}{{\spec K[t]/(t^{m+1})}}
\newcommand{\sT}{{\spec K[[t]]}}
\newcommand{\tm}{{k[t]/(t^{m+1})}}
\newcommand{\D}{{\Delta}}
\newcommand{\fa}{{\frak{a}}}
\newcommand{\fb}{{\frak{b}}}
\newcommand{\fm}{{\frak{m}}}
\newcommand{\fp}{{\frak{p}}}
\newcommand{\bx}{{\bold{x}}}
\newcommand{\bz}{{\bold{z}}}
\newcommand{\oa}{{\overline{\overline \fa}}}
\newcommand{\ta}{{\tilde{ \fa}}}
\newcommand{\tb}{{\tilde{ \fb}}}
\newcommand{\ia}{{\overline\fa}}
\newcommand{\ib}{{\overline\fb}}
\newcommand{\ch}{{\operatorname{char}}}

\newcommand{\ob}{{\overline{\overline\fb}}}
\let \cedilla =\c
\renewcommand{\c}[0]{{\mathbb C}}  
\let \crossedo =\o
\renewcommand{\o}[0]{{\mathcal O}} 
\newcommand{\z}[0]{{\mathbb Z}}
\newcommand{\n}[0]{{\mathbb N}}

\let \ringaccent=\r  
\renewcommand{\r}[0]{{\mathbb R}} 

\renewcommand{\a}[0]{{\mathbb A}} 

\newcommand{\h}[0]{{\mathbb H}}
\newcommand{\p}[0]{{\mathbb P}}
\newcommand{\f}[0]{{\mathbb F}}
\newcommand{\q}[0]{{\mathbb Q}}
\newcommand{\map}[0]{\dasharrow}
\newcommand{\qtq}[1]{\quad\mbox{#1}\quad}
\newcommand{\sing}[0]{\operatorname{Sing}}

\newcommand{\cont}{{\operatorname{Cont}}}

\newcommand{\onto}[0]{\twoheadrightarrow}

\def\into{\DOTSB\lhook\joinrel\rightarrow}
\def\to {\longrightarrow}

\newtheorem{thm}{Theorem}

\newtheorem{lem}[thm]{Lemma}
\newtheorem{cor}[thm]{Corollary}
\newtheorem{prop}[thm]{Proposition}
\newtheorem{problem}[thm]{Problem}

\theoremstyle{definition}
\newtheorem{defn}[thm]{Definition}

\newtheorem{say}[thm]{}
\newtheorem{exmp}[thm]{Example}
\newtheorem{conj}[thm]{Conjecture}

\newtheorem{rem}[thm]{Remark}

\theoremstyle{remark}
\newtheorem{case}{Case}

\markboth{\hfill S. Ishii\hfill}{\hfill
Supplement to Fujino's paper  \hfill}

\begin{abstract}
Let $E$ be the essential part of the exceptional locus of a good resolution of an isolated, log canonical singularity of index one. We describe the dimension of the dual complex of $E$ in terms of the Hodge 
type of $H^{n-1}(E,\o_E)$, which is one of the main results of the paper \cite{f} of Fujino.
Our proof uses only an elementary classical method, while Fujino's argument depends on the recent development in minimal model theory.
\end{abstract}

\noindent
In this paper, a normal singularity
 $(X,x)$ of dimension $n\geq 2$ is always assumed to be  isolated, strictly log canonical of index one, where a strictly log canonical singularity means a log canonical and not log terminal singularity.
Let $f:Y\to X$ be a good resolution ({\it i.e}., a resolution with the simple normal crossing exceptional divisor) of the singularity of $X$. We have
$$K_Y=f^*K_X- E+F,$$
where $E=E_{red}>0$ and $F\geq 0$ have no common components.
The divisor $E$ is called the essential part of the exceptional divisor on the resolution.
For a simple normal crossing divisor $E$, we associate a simplicial complex $\Gamma_E$
called the dual complex in a canonical way.
Fujino defines an invariant $\mu(X,x)$ and it turns out to be
$$\mu=\mu(X,x)=\min\{\dim W \mid W \mbox{is\ a\ stratum\ of\ } E\}$$
(see \cite[4.11]{f}).
Note that $\dim \Gamma_E=n-\mu-1$.

On the other hand,
we define the Hodge type of the singularity $(X,x)$ in the following way:
Since
$$\c=H^{n-1}(E, \o_E)\simeq Gr_F^0H^{n-1}(E,\c)\simeq \bigoplus_{i=0}^{n-1}H_{n-1}^{0,i}(E),$$
there is unique $i$ such that $H_{n-1}^{0,i}(E)\neq 0$, where $H_{n-1}^{0,i}(E)$ is $(0,i)$-Hodge component of $H^{n-1}(E,\c)$ and $F$ is the Hodge filtration.
In this case, we call the singularity $(X,x)$ of type $(0,i)$.
We can easily prove that the type is independent of the choice of resolutions (\cite{i1}).

One of the main results (Theorem 5.5) in \cite{f}  states that for $(X,x)$ of type $(0,i)$, the equality $\mu(X,x)=i$ holds.
Theorem \ref{th} below states the same conclusion and its
 proof  was privately communicated by the author to Fujino in 1999 (cited as  [I3] in the reference list of \cite{f}).
The author thinks that it is reasonable to publish the original proof   as a supplement to Fujino's article, because her original proof is simpler and used only classical method, while Fujino uses recent results in minimal model theory.

\begin{thm}
\label{th}
 {\sl  Let \( E \) be the essential part of the exceptional divisor of a good resolution \( Y\to X \)
 of an \( n \)-dimensional 
 isolated strictly log canonical singularity $(X,x)$.%
  If the Hodge type is of \( (0,i) \), then
  \( \dim \Gamma_{E}= n-i-1 \).}
\end{thm}  

The following lemma appeared in \cite[Lemma 7.4.9]{i2}.
As it is written in Japanese, we write the proof down here for the non-Japanese readers.
\begin{lem}
\label{l}
  Let $E$ be a simple normal crossing divisor on an $n$-dimensional non-singular variety.
  If $H_{n-1}^{0,i}(E)\neq 0$, then $\dim \Gamma_E\geq n-i-1$.
  \end{lem}

\begin{pf} After renumbering the suffixes if necessary, we prove that  there exist $n-i$ irreducible components 
$E_1,\ldots, E_{n-i}$ such that $E_1\cap\cdots\cap E_{n-i}\neq \emptyset$.
Let $E'$ be a minimal subdivisor of $E$ such that  $H_{n-1}^{0,i}(E')
\neq 0$. 
If $E'$ is irreducible, then it is a non-singular variety of dimension $n-1$, therefore we obtain
$i=n-1$ by the basic fact in mixed Hodge theory (see for example \cite[Theorem 7.1.6]{i2}).
Therefore, 
$$\dim \Gamma_E\geq 0=n-(n-1)-1,$$
{\it i.e}., the required inequality becomes trivial. 
If $E'$ is not irreducible, take an irreducible component $E_1< E'$ and decompose
$E'$ as  $E'=
E_1+E^\vee_1$.
Then by the minimality of $E'$, we have $H_{n-1}^{0,i}(E_1)=H_{n-1}^{0,i}(E_1^\vee)=0$.
Consider the exact sequence:
$$H^{n-2}(E_1\cap E_1^\vee,\c)\to H^{n-1}(E',\c)\to H^{n-1}(E_1,\c)\oplus H^{n-1}(E_1^\vee,\c).$$
By the above vanishing, the $(0,i)$-component of the center term comes from the left term,
therefore $i\leq n-2$ and  $H_{n-2}^{0,i}(E_1\cap E_1^\vee)\neq 0$.
  
  Take $E_1^\dagger$, a minimal subdivisor of $E_1^\vee$ such that   
$H_{n-2}^{0,i}(E_1\cap E_1^\dagger)\neq 0$.
If $E_1\cap E_1^\dagger$ is irreducible, then it is a non-singular variety of dimension $n-2$,
therefore we obtain $i=n-2$ by the basic fact in mixed Hodge theory.
Therefore, 
$$\dim \Gamma_E\geq \dim \Gamma_{E_1+E_1^\dagger}\geq 1=n-(n-2)-1,$$
{\it i.e}., the required inequlatiy holds.
If $E^\dagger=E_1\cap E_1^\dagger$ is not irreducible, take an irreducible component $E_2$ of $E^\dagger$ such that the decomposition $E^\dagger=E_2+E_2^\vee$ gives a non-trivial 
decomposition $E_1\cap E_1^\dagger=E_1\cap E_2+E_1\cap E_2^\vee$.
By the same argument as above, we obtain $i\leq n-3$ and  $H_{n-3}^{0,i}(E_1\cap E_2\cap E_2^\vee)\neq
0$.
Continue this procedure successively 
until we eventually obtain 
$$H_i^{0,i}(E_1\cap E_2\cap\cdots \cap E_{n-i-1}\cap E_{n-i-1}^\vee)\neq 0,$$
which yields 
$E_1\cap E_2\cap\cdots \cap E_{n-i-1}\cap E_{n-i-1}^\vee\neq \emptyset$.
\end{pf}

\vskip.5truecm
\noindent   
{\sl Proof of Theorem \ref{th}}. 
 The inequality \( \geq \) is proved in Lemma \ref{l}.
  Assume the strict inequality.
  Then there exist components \( E_{1},\ldots \ E_{s} \), \( (s> 
  n-i) \) such that \( C:=E_{1}\cap\ldots\cap E_{s}\neq \emptyset \).
  We may assume that \( E_{j}\cap C =\emptyset \) for any \( E_{j} \)
  \( (j>s) \).
  Let \( \varphi:Y'\to Y \) be the blow-up at \( C \),
  \( E' \)  the reduced total pull-back of \( E \),
  \( E_{0} \) the exceptional divisor for \( \varphi \)
  and \( E'_{j} \) the proper transform of \( E_{j} \).
  Then \( E' \) is again the essential part on $Y'$ and \( E' \) itself is 
  a minimal subdivisor of \( E' \) such that \( H^{0,i}_{n-1}(E')\neq 0 \) by \cite[Corollary 3.9]{i1}.
  Make the 
  procedure of the proof of the lemma with taking \( E_{0} \)
  as \( E_{1} \) in the lemma.
  Then we obtain \( E'_{1},\ldots,E'_{n-i-1} \) (by renumbering the 
  suffices \( 1,\ldots, s \)) such that
  \( H^{0,i}_{i}(E_{0}\cap E'_{1}\cap\ldots\cap E'_{n-i-1})\neq 0 \).
  On the other hand, \( i \)-dimensional variety \( E_{0}\cap E'_{1}\cap\ldots\cap E'_{n-i-1} \)
  is a \( {\Bbb P}^{s-n+i} \)-bundle over \( C \), because  it is  the exceptional divisor of the blow up of an
  $(i+1)$-dimensional variety $E_1\cap \cdots\cap E_{n-i-1}$ with the $(n-s)$-dimensional center $C$. By the assumption on $s$, we note that $s-n+i>0$.
  Hence we have $H^i(E_{0}\cap E'_{1}\cap\ldots\cap E'_{n-i-1}, \o)=0$. In particular
  $$H_i^{0, i}(E_{0}\cap E'_{1}\cap\ldots\cap E'_{n-i-1})=0,$$
  a contradiction.
$\square$
%

 \end{document}